\RequirePackage{ifpdf}
\ifpdf 
\documentclass[pdftex]{sigma}
\else
\documentclass{sigma}
\fi

\numberwithin{equation}{section}

\newtheorem{thm}{Theorem}[section]
\newtheorem{cor}[thm]{Corollary}
\newtheorem{lem}[thm]{Lemma}
\newtheorem{prop}[thm]{Proposition}
{\theoremstyle{definition}
\newtheorem{rem}[thm]{Remark}
\newtheorem{defn}[thm]{Definition}}

\newcommand{\C}{\mathbb C}
\newcommand{\CP}{\mathbb{CP}}

\newcommand{\Z}{\mathbb Z}
\def\Gamm#1{\Gamma_{\!\!#1}}
\def\Gampq{\Gamm{p,q}}
\newcommand{\II}{\mathord{I\!I}}
\def\disj{\uplus}
\DeclareMathOperator{\GL}{GL}
\DeclareMathOperator{\PGL}{PGL}

\def\rowvec#1#2{\begin{matrix}\smash{\bigl(}#1&#2\smash{\bigr)}\\\mathstrut&\end{matrix}}

\begin{document}

\allowdisplaybreaks

\renewcommand{\thefootnote}{$\star$}

\renewcommand{\PaperNumber}{088}

\FirstPageHeading

\ArticleName{An Isomonodromy Interpretation\\ of the   Hypergeometric Solution
  of the Elliptic \\ Painlev\'e Equation (and Generalizations)\footnote{This paper is a
contribution to the Special Issue ``Relationship of Orthogonal Polynomials and Special Functions with Quantum Groups and Integrable Systems''. The
full collection is available at
\href{http://www.emis.de/journals/SIGMA/OPSF.html}{http://www.emis.de/journals/SIGMA/OPSF.html}}}

\ShortArticleName{Isomonodromy Interpretation of Hypergeometric Solution of Elliptic Painlev\'e}

\Author{Eric M.~RAINS}

\AuthorNameForHeading{E.M.~Rains}

\Address{Department of Mathematics, California Institute of   Technology,\\
 1200 E.~California Boulevard, Pasadena, CA 91125, USA}
\Email{\href{mailto:rains@caltech.edu}{rains@caltech.edu}}

\ArticleDates{Received April 25, 2011, in f\/inal form September 06, 2011;  Published online September 09, 2011}

\Abstract{We construct a family of second-order linear dif\/ference equations
parametrized by the hypergeometric solution of the elliptic Painlev\'e
equation (or higher-order analogues), and admitting a large family of
monodromy-preserving deformations.  The solutions are certain semiclassical
biorthogonal functions (and their Cauchy transforms), biorthogonal with
respect to higher-order analogues of Spiridonov's elliptic beta integral.}


\Classification{33E17; 34M55; 39A13}

\section{Introduction}

In \cite{SakaiH:2001}, Sakai introduced an elliptic analogue of the
Painlev\'e equations, including all of the known discrete (and continuous)
Painlev\'e equations as special cases.  Unfortunately, although Sakai's
construction is quite natural and geometric, it does not ref\/lect the most
important role of the ordinary Painlev\'e transcendents, namely as parameters
controlling monodromy-preserving deformations.

As with the ordinary Painlev\'e equations, the elliptic Painlev\'e equation
admits a special class of ``hypergeometric'' solutions
\cite{KajiwaraK/MasudaT/NoumiM/OhtaY/YamadaY:2003,recur} that in the most
general case can be expressed via $n$-dimensional contour integrals with
integrands expressed in terms of elliptic Gamma functions.  It is thus
natural, as a f\/irst step in constructing an isomonodromy interpretation of
the elliptic Painlev\'e equation, to attempt to understand that
interpretation in the hypergeometric case (and thus gain insight to the
general case).  Note that we want to understand the hypergeometric case for
all $n\ge 1$, to avoid the possibility that the small $n$ cases might
dif\/fer from the general Painlev\'e case in some qualitatively signif\/icant
way.  (For instance, \cite{KitaevAV:2000} considers the isomonodromy
interpretation of the usual ${}_2F_1$ (corresponding to $n=1$ in our
setting), but this is simplif\/ied greatly from the general Painlev\'e VI
case by the fact that not only the monodromy but the equation itself can be
taken to be triangular.)

In the present work, we do precisely that: associated to each elliptic
hypergeometric solution of the elliptic Painlev\'e equation, we construct a
corresponding second-order linear dif\/ference equation that admits a family
of discrete ``monodromy-preserving'' deformations.  (In fact, the
construction works equally well for higher-order analogues of the relevant
elliptic hypergeometric integrals, which should correspond to special
solutions of ``elliptic Garnier equations''.) The construction is based on
an analogue of the approach in
\cite{MagnusAP:1995,ForresterPJ/WitteNS:2006}.  There, a linear
dif\/ferential equation deformed by the hypergeometric case of the Painlev\'e
VI equation is constructed as a~dif\/ferential equation satisf\/ied by a family
of ``semiclassical'' (bi-)orthogonal polynomials.  Our construction is much
the same, although there are several technical issues to overcome.

The f\/irst such issue is, simply put, to understand precisely what it means
for a deformation of an elliptic dif\/ference equation to preserve monodromy,
or even what the monodromy of an elliptic dif\/ference equation is.  While we
give only a partial answer to this question, we do def\/ine (in Section~\ref{section2}
below; note that many of the considerations there turn out to have been
anticipated by Etingof in~\cite{EtingofPI:1993}) a weakened form of
monodromy that, while somewhat weaker than the analogous notions at the
$q$-dif\/ference~\cite{EtingofPI:1995} and lower~\cite{BirkhoffGD:1911}
levels, is still strong enough to give a reasonably rigid notion of
isomonodromy deformation.  Indeed, two elliptic dif\/ference equations have
the same weak monodromy if\/f the corresponding dif\/ference modules (see~\cite{vanderPutM/SingerMF:1997}) are isomorphic; the same holds for
ordinary dif\/ference equations, even relative to the stronger notion of
monodromy \cite{BorodinA:2004}.  The key observation is that a fundamental
matrix for a $p$-elliptic $q$-dif\/ference equation is {\em also} a~fundamental matrix for a $q$-elliptic $p$-dif\/ference equation; this latter
equation (up to a certain equivalence relation) plays the role of the
monodromy.  (The result is similar to the notion of monodromy introduced by
Krichever in~\cite{KricheverIM:2004}; while our notion is weakened by an
equivalence relation, it avoids any assumptions of genericity.)

In Section~\ref{section3}, we develop the theory of semiclassical elliptic biorthogonal
functions, functions biorthogonal with respect to a density generalizing
Spiridonov's elliptic beta integral~\cite{SpiridonovVP:2001a} by adding $m$
additional pairs of parameters.  The key observation is that such functions
can be constructed as higher-order elliptic Selberg integrals of a special
form; in addition, their ``Cauchy transforms'' can also be so written.
This gives rise to several nice relations between these functions, which we
describe.  Most important for our purposes is their behavior under
$p$-shifts; the biorthogonal functions themselves are $p$-elliptic, but if
we include the Cauchy transforms, the overall action is triangular.  We can
thus construct from these functions a $2\times 2$ matrix which satisf\/ies a
triangular $q$-elliptic $p$-dif\/ference equation, analogous to the
Riemann--Hilbert problem associated to orthogonal polynomials (\cite[\S~3.4]{FokasAS/ItsAR/KitaevAV:1992}; see also~\cite{DeiftPA:1999} for a
general exposition).  By the theory of Section~\ref{section2}, this immediately gives
rise to a $p$-elliptic $q$-dif\/ference equation, and symmetries of the
$p$-dif\/ference equation induce monodromy-preserving deformations of the
$q$-dif\/ference equation.

Finally, in Section~\ref{section4}, we compute this dif\/ference equation and the
associated deformations.  Although we cannot give a closed form expression
for the dif\/ference equation, we are able at least to determine precisely
where the dif\/ference equation is singular, and at each such point, compute
the value (or residue, as appropriate) of the shift matrix.  Together with
the fact that the coef\/f\/icients are meromorphic $p$-theta functions, this
data suf\/f\/ices to (over)determine the shift matrix.

In a followup paper~\cite{partII}, with Arinkin and Borodin, we will
complete the isomonodromy interpretation of the elliptic Painlev\'e
equation by applying the ideas of~\cite{ArinkinD/BorodinA:2006} to show
that any dif\/ference equation having the same structure as the ones
constructed below admits a corresponding fa\-mi\-ly of monodromy-preserving
deformations, and moreover that (when $m=1$) Sakai's rational surface can
be recovered as a moduli space of such dif\/ference equations.  A rather
dif\/ferent geo\-metric approach to such an interpretation (via a Lax pair) for
the $m=1$ case has been given in~\cite{YamadaY:2009}.

\subsection*{Notation}

The {\em elliptic Gamma function} \cite{RuijsenaarsSNM:1999} is def\/ined for
complex numbers $p$, $q$, $z$ with $|p|,|q|<1$, $z\ne 0$, by
\[
\Gampq(z) := \prod_{0\le i,j} \frac{1-p^{i+1}q^{j+1}/z}{1-p^i q^j z},
\]
and satisf\/ies the ref\/lection relation
\[
\Gampq(pq/z) = \Gampq(z)^{-1}
\]
as well as the shift relations
\begin{gather*}
\Gampq(pz)  = \theta_q(z)\Gampq(z),\qquad
\Gampq(qz)  = \theta_p(z)\Gampq(z),
\end{gather*}
where the function
\[
\theta_p(z) := \prod_{0\le i} \big(1-p^{i+1}/z\big)\big(1-p^i z\big)
\]
satisf\/ies
\[
\theta_p(z) = -z \theta_p(1/z) = \theta_p(p/z),
\]
so that
\[
\Gampq(pqz)\Gampq(z) = -z^{-1}\Gampq(pz)\Gampq(qz).
\]
By convention, multiple arguments to a Gamma or theta function represent a
product; thus, for instance
\[
\Gampq(u_0 z^{\pm 1}) = \Gampq(u_0 z)\Gampq(u_0/z).
\]
We will also make brief use of the third-order elliptic Gamma function
\[
\Gamma^+_{p,q,t}(x)
:=
\prod_{0\le i,j,k} \big(1-p^iq^jt^k x\big)\big(1-p^{i+1}q^{j+1}t^{k+1}/x\big),
\]
which satisf\/ies
\begin{gather*}
\Gamma^+_{p,q,t}(tx) =\Gampq(x)\Gamma^+_{p,q,t}(x),\qquad
\Gamma^+_{p,q,t}(pqt/x)  = \Gamma^+_{p,q,t}(x);
\end{gather*}
for our purposes, this appears only as a normalization factor relating the
order~1 elliptic Selberg integral to the hypergeometric tau function for
elliptic Painlev\'e.

\section{Elliptic dif\/ference equations}\label{section2}

Let $p$ be a complex number with $|p|<1$.  A (meromorphic) $p$-theta
function of multiplier $\alpha z^k$ is a meromorphic function $f(z)$ on
$\C^*:=\C\setminus \{0\}$ with the periodicity property $f(pz)=\alpha z^k
f(z)$.  (To justify this def\/inition, observe that the composition
$f(\exp(2\pi\sqrt{-1} t))$ is meromorphic on~$\C$, periodic with period
$1$, and quasi-periodic with period $\log(p)/2\pi\sqrt{-1}$; in other
words, it is a theta function in the usual sense.)  The canonical example
of such a function is $\theta_p(z)$, a~holomorphic $p$-theta function with
multiplier $-z^{-1}$; indeed, any holomorphic $p$-theta function can be
written as a product of functions~$\theta_p(uz)$, and any meromorphic
$p$-theta function as a~ratio of such products.  In the special case of
multiplier~1, the function is called $p$-elliptic, for similar reasons.  By
standard convention, a $p$-theta function, if not explicitly allowed to be
meromorphic, is holomorphic; however, $p$-elliptic functions are always
allowed to be meromorphic (since a holomorphic $p$-elliptic function is
constant).

Let $q$ be another complex number with $|q|<1$, such that $p^\Z\cap
q^\Z=\varnothing$.

\begin{defn}
A {\em $p$-theta $q$-difference equation of multiplier $\mu(z)=\alpha z^k$} is
an equation of the form
\[
v(qz) = A(z) v(z),
\]
where $A(z)$ is a {\em nonsingular meromorphic matrix} (a square matrix,
each coef\/f\/icient of which is meromorphic on $\C^*$, and the determinant of
which is not identically 0), called the {\em shift matrix} of the equation,
such that
\[
A(pz)=\mu(z)A(z),
\]
so in particular the coef\/f\/icients of $A$ are meromorphic $p$-theta functions of
multiplier $\mu(z)$.  Si\-milarly, a $p$-elliptic $q$-dif\/ference equation is
a $p$-theta $q$-dif\/ference equation of multiplier $1$.
\end{defn}

We will refer to the dimension of the matrix $A$ as the {\em order} of the
corresponding dif\/ference equation.  We note the following fact about
nonsingular meromorphic matrices.

\begin{prop}
Let $M(z)$ be a nonsingular meromorphic matrix.  Then $M(z)^{-1}$ is also a~nonsingular meromorphic matrix, and if the coefficients of $M(z)$ are
meromorphic $p$-theta functions of multiplier $\mu(z)$, then those of
$M(z)^{-1}$ are meromorphic $p$-theta functions of multiplier~$\mu(z)^{-1}$.
\end{prop}

\begin{proof}
Indeed, the coef\/f\/icients of the adjoint matrix $\det(M(z))M(z)^{-1}$ are
minors of $M(z)$, and thus, as polynomials in meromorphic functions, are
meromorphic; this continues to hold after multiplying by the meromorphic
function $\det(M(z))^{-1}$.  For the second claim, if
\[
M(pz)=\mu(z)M(z),
\]
then
\begin{equation*}
M(pz)^{-1} = \mu(z)^{-1}M(z)^{-1}.\tag*{\qed}
\end{equation*}
\renewcommand{\qed}{}
\end{proof}

\begin{defn}
Let $v(qz)=A(z)v(z)$ be a $p$-theta $q$-dif\/ference equation.  A {\em
  meromorphic fundamental matrix} for this equation is a~nonsingular
meromorphic matrix $M(z)$ satisfying
\[
M(qz) = A(z)M(z).
\]
\end{defn}

It follows from a theorem of Praagman \cite[Theorem~3]{PraagmanC:1986} that
for any nonsingular meromorphic matrix $A(z)$, there exists a nonsingular
meromorphic matrix $M(z)$ satisfying $M(qz)=A(z)M(z)$ (this is the special
case of the theorem in which the discontinuous group acting on $\CP^1$ is
that generated by multiplication by $q$).  In particular, any $p$-theta
$q$-dif\/ference equation admits a meromorphic fundamental matrix.  In the
case of a f\/irst order equation, we can explicitly construct such a matrix.

\begin{prop}
Any first order $p$-theta $q$-difference equation admits a meromorphic
fundamental matrix.
\end{prop}

\begin{proof}
For any nonzero meromorphic $p$-theta function $a(z)$, we need to construct
a nonzero meromorphic function $f(z)$ such that
\[
f(qz) = a(z) f(z).
\]
Since $a(z)$ can be factored into functions $\theta_p(uz)$, it suf\/f\/ices to
consider the case $a(z)=\theta_p(uz)$, with meromorphic solution
\[
f(z) = \Gampq(uz);
\]
this includes the case $a(z)=b z^k$ by writing
\begin{equation*}
b z^k =
\frac{\theta_p(-bz)\theta_p(-pz)^{k-1}}{\theta_p(-bpz)\theta_p(-z)^{k-1}}.
\tag*{\qed}
\end{equation*}
\renewcommand{\qed}{}
\end{proof}

We note in particular that, since the elliptic Gamma function is
symmetrical in $p$ and $q$, the solution thus obtained for a f\/irst order
$p$-theta $q$-dif\/ference equation also satisf\/ies a $q$-theta $p$-dif\/ference
equation.  This is quite typical, and in fact we have the following result.

\begin{lem}\label{lem:q_from_p}
Let $v(qz)=A(z)v(z)$ be a $p$-theta $q$-difference equation of multiplier
$\mu(z)$, and let $M(z)$ be a meromorphic fundamental matrix for this
equation.  Then there exists a $($unique$)$ $q$-theta $p$-difference equation
of multiplier $\mu(z)$ for which $M(z)^t$ is a fundamental matrix.
\end{lem}

\begin{proof}
An equation $w(pz)=C(z)w(z)$ with fundamental matrix $M(z)^t$ satisf\/ies
\[
M(pz)^t = C(z) M(z)^t,
\]
and thus, since $M(z)$ is nonsingular, we can compute
\[
C(z) = M(pz)^t M(z)^{-t}.
\]
(Here $M^{-t}$ denotes the inverse of the transpose of $M$.)  This matrix is
meromorphic, and satisf\/ies
\begin{equation*}
C(qz) = M(pqz)^t M(qz)^{-t}
      = M(pz) A(qz)^t A(z)^{-t} M(z)^{-t}
      = \mu(z) C(z).
\tag*{\qed}
\end{equation*}
\renewcommand{\qed}{}
\end{proof}

By symmetry, we obtain the following result.

\begin{thm}
Let $M(z)$ be a nonsingular meromorphic matrix.  Then the following are
equivalent:
\begin{itemize}\itemsep=0pt
\item[$(1)$] $M(z)$ is a meromorphic fundamental matrix for some $p$-theta
  $q$-difference equation;
\item[$(2)$] $M(z)^t$ is a meromorphic fundamental matrix for some $q$-theta
  $p$-difference equation;
\item[$(1')$] $M(z)^{-t}$ is a meromorphic fundamental matrix for some $p$-theta
  $q$-difference equation;
\item[$(2')$] $M(z)^{-1}$ is a meromorphic fundamental matrix for some $q$-theta
  $p$-difference equation,
\end{itemize}
as are the corresponding statements with ``some'' replaced by ``a unique''.
Furthermore, if the above conditions hold, the multipliers of the
difference equations of~$(1)$ and~$(2)$ agree, and are inverse to those of~$(1')$
and~$(2')$.
\end{thm}

\begin{rem}
In the elliptic case, the above observations were made by Etingof
\cite{EtingofPI:1993}, who also noted that the associated $q$-elliptic
$p$-dif\/ference equation can be thought of as the monodromy of $M$.
\end{rem}

Given a $p$-theta $q$-dif\/ference equation, the corresponding meromorphic
fundamental matrix is by no means unique, and thus we obtain a whole family
of related $q$-theta $p$-dif\/ference equations.  There is, however, a
natural equivalence relation on $q$-theta $p$-dif\/ference equations such
that any $p$-theta $q$-dif\/ference equation gives rise to a well-def\/ined
equivalence class.  First, we need to understand the extent to which the
fundamental matrix fails to be unique.

\begin{lem}
Let $M(z)$ and $M'(z)$ be fundamental matrices for the same $p$-theta
$q$-difference equation $v(qz)=A(z)v(z)$.  Then
\[
M'(z) = M(z) D(z)^t
\]
for some nonsingular meromorphic matrix $D(z)$ with $q$-elliptic coefficients.
\end{lem}

\begin{proof}
Certainly, there is a unique meromorphic matrix $D(z)$ with
$M'(z)=M(z)D(z)^t$, and comparing determinants shows it to be nonsingular.
It thus remains to show that $D(z)$ has $q$-elliptic coef\/f\/icients, or
equivalently that $D(qz)=D(z)$.  As in the proof of Lemma~\ref{lem:q_from_p}, we can write
\[
A(z) = M(qz)M(z)^{-1} = M'(qz)M'(z)^{-1},
\]
and thus
\[
D(qz)^t D(z)^{-t} = M(qz)^{-1} M'(qz) M'(z)^{-1} M(z) = 1,
\]
as required.
\end{proof}

\begin{thm}\label{thm:weak_monodromy_defined}
Define an equivalence relation on $q$-theta $p$-difference equations by
saying
\[
\bigl[ v(pz)=C(z)v(z) \bigr] \cong \bigl[ v(pz) = C'(z) v(z)\bigr]
\]
iff there exists a nonsingular $q$-elliptic matrix $D(z)$ such that
\[
C'(z) D(z) = D(pz) C(z).
\]
Then the set of $q$-theta $p$-difference equations associated to a given
$p$-theta $q$-difference equation is an equivalence class.
\end{thm}

\begin{proof}
Let $M(z)$ be a meromorphic fundamental matrix for the $p$-theta
$q$-dif\/ference equation $v(qz)=A(z)v(z)$, and associated $q$-theta
$p$-dif\/ference equation $w(pz)=C(z)w(z)$.  If $M'(z)$ is another
meromorphic fundamental matrix for the $q$-dif\/ference equation, with
associated $p$-dif\/ference equation $w(pz)=C'(z)w(z)$, then
\[
M'(z) = M(z)D(z)^t,
\]
and thus
\[
C'(z)
=
M'(pz)^t M'(z)^{-t}
=
D(pz) M(pz)^t M(z)^{-t} D(z)^{-1}
=
D(pz) C(z) D(z)^{-1}.
\]
Conversely, if
\[
C'(z) D(z) = D(pz)C(z),
\]
then $M(z)D(z)^t$ is a fundamental matrix for a $q$-dif\/ference equation
with associated $p$-dif\/ference equation $w(pz)=C'(z)w(z)$.
\end{proof}

\begin{defn}
The {\em weak monodromy} of a $p$-theta $q$-dif\/ference equation is the
associated equiva\-lence class of $q$-theta $p$-dif\/ference equations.  Two
$p$-theta $q$-dif\/ference equations are {\em iso\-monodromic} if they have
the same weak monodromy.
\end{defn}

\begin{thm}
The $p$-theta $q$-difference equations $v(qz)=A(z)v(z)$, $v(qz)=A'(z)v(z)$
are isomonodromic iff there exists a nonsingular $p$-elliptic matrix $B(z)$
such that
\[
A'(z) B(z) = B(qz) A(z).
\]
\end{thm}

\begin{proof}
Choose a $q$-theta $p$-dif\/ference equation $w(qz)=C(z)w(z)$ representing
the weak mo\-no\-dromy of the f\/irst equation.  The equations are
isomonodromic if\/f $w(qz)=C(z)w(z)$ represents the weak monodromy of the
second equation, if\/f the two equations have fundamental matrices satisfying
$M(qz)=M(z)C(z)$.  But by Theorem \ref{thm:weak_monodromy_defined}
(swapping $p$ and $q$), this holds if\/f $A'(z)B(z)=B(qz)A(z)$ for some
$p$-elliptic matrix $B(z)$.
\end{proof}

\begin{rem}
Compare \cite{BorodinA:2004}, where the analogous result is proved for
dif\/ference equations, relative to Birkhof\/f's \cite{BirkhoffGD:1911} notion of
monodromy.
\end{rem}

\begin{cor}
The map from isomonodromy classes of $p$-theta $q$-difference equations to
their weak monodromies is well-defined, and inverse to the map from
isomonodromy classes of $q$-theta $p$-difference equations to {\em their}
weak monodromies.
\end{cor}

\begin{rem}
The isomonodromy equivalence relation is also quite natural from
the perspective of the general theory of dif\/ference equations (see, e.g.,
\cite{vanderPutM/SingerMF:1997}); to be precise, two $p$-theta
$q$-dif\/ference equations are isomonodromic if\/f they induce isomorphic
dif\/ference modules.  The latter fact induces a natural isomorphism between
their dif\/ference Galois groups (at least when the latter are def\/ined, i.e.,
when the equations are elliptic), as can be seen directly from the
interpretation of dif\/ference Galois groups via Tannakian categories.  This
preservation of Galois groups seems to be what is truly intended by the
word ``isomonodromy'', even in the dif\/ferential setting.  For instance, for
non-Fuchsian equations, where the monodromy group conveys relatively little
information, one only obtains the relevant Painlev\'e equations by
insisting that the corresponding deformations should preserve Stokes data
as well.
\end{rem}

It will be convenient in the sequel to introduce a slightly weaker
equivalence relation.

\begin{defn}
Two $p$-theta $q$-dif\/ference equations are {\em theta-isomonodromic} if
there exists a~nonsingular meromorphic $p$-theta matrix $B(z)$ such that
the shift matrices $A(z)$, $A'(z)$ of the equations satisfy $A'(z) B(z) =
B(qz) A(z)$.
\end{defn}

\begin{thm}
Two $p$-theta $q$-difference equations are theta-isomonodromic iff their
weak mo\-nodromies agree up to multiplication of the shift matrix by a factor
of the form $a z^k$.
\end{thm}

Note that $v(qz)=A(z)v(z)$ and $v(qz)=A(qz)v(z)$ are theta-isomonodromic
with $B(z)=A(z)$.
\medskip

\begin{rem}
Though this equivalence relation no longer preserves Galois groups,
even if both equations are elliptic, it comes quite close to doing
so.  Indeed, the Galois group of an $n$-th order equation is naturally
a subgroup of $\GL_n$, and we may thus consider its image in~$\PGL_n$,
which one might call the {\em projective} Galois group.  Since
$\PGL_n$ is the image of $\GL_n$ under the adjoint representation,
one f\/inds that the projective Galois group of the equation with
shift matrix $A(z)$ can be identif\/ied with the ordinary Galois group
of the equation with shift matrix $A(z)\otimes A(z)^{-t}$.  If two
equations are theta-isomonodromic, their images under the adjoint
representation are thus fully isomonodromic, and thus the original
equations had the same projective Galois groups.  This even extends
to $p$-theta $q$-dif\/ference equations once we observe that the image
of such an equation under the adjoint representation is elliptic,
and thus the projective Galois groups of such equations are still
well-def\/ined.
\end{rem}

\begin{rem}
The relation of isomonodromy to Galois groups suggests some further
questions, which are in the main outside the scope of the current paper,
but seem to merit a brief mention nonetheless.

First, since (theta-)isomonodromic equations have isomorphic (projective)
Galois groups, it is natural to ask whether one can recover the
(projective) Galois group more directly from the weak monodromy.  Since the
weak monodromy is itself a isomonodromy class, any two representatives of
the weak monodromy have the same Galois group, and one would expect that
group to be related to the original Galois group.  It can be shown
(Etingof, personal communication) that in fact the groups are naturally
isomorphic, with dual associated representations.  Thus, for instance, the
fact that the dif\/ference equations we will be considering have triangular
weak monodromy implies that they have solvable Galois group.  (This also
follows immediately from the fact that, by construction, they have theta
function solutions.)

Another natural question is whether there exists a stronger notion of
monodromy; for rational $q$-dif\/ference equations with suf\/f\/iciently nice
singularities, there is a well-def\/ined notion of monodromy, an associated
nonsingular $q$-elliptic matrix the nonsingular values of which generate a
Zariski dense subgroup of the Galois group~(\cite{EtingofPI:1995}; see also
Chapter~12 of~\cite{vanderPutM/SingerMF:1997}).  Krichever~\cite{KricheverIM:2004} def\/ines an analogous matrix for generic dif\/ference
equations with theta function coef\/f\/icients (although the relation to the
Galois group is again unclear); although Krichever's genericity assumptions
explicitly exclude the situation we consider above (raising the question of
whether there is an analogue in our setting), his monodromy is again a
dif\/ference equation with theta function coef\/f\/icients.  This suggests that
the rational $q$-dif\/ference notion of monodromy should correspond at the
elliptic level to a representative of our weak monodromy, and thus suggests
the question of whether given a $p$-elliptic $q$-dif\/ference equation, there
exists a representative of its weak monodromy such that the nonsingular
values of the corresponding $C$ matrix are Zariski dense in its Galois
group.
\end{rem}

\section{Semiclassical biorthogonal elliptic functions}\label{section3}

In \cite{SpiridonovVP:2001b}, Spiridonov constructed a family of elliptic
hypergeometric functions biorthogonal with respect to the density of the
elliptic beta integral:
\[
\frac{(p;p)(q;q)}{2}
\int_C \frac{\prod\limits_{0\le r<6} \Gampq(u_r z^{\pm 1})}{\Gampq(z^{\pm 2})}
\frac{dz}{2\pi\sqrt{-1}z}
=
\prod_{0\le r<s<6} \Gampq(u_r u_s),
\]
where the parameters satisfy the {\em balancing condition}
\[
\prod_{0\le r<6} u_r = pq,
\]
and the (possibly disconnected, but closed) contour is chosen to be
symmetrical under $z\mapsto 1/z$, and to contain all points of the form
$p^i q^j u_r$, $i,j\ge 0$, $0\le r<6$, or more precisely, all poles of the
integrand of that form.

If we view this as the ``classical'' case, then this suggests, by analogy
with~\cite{MagnusAP:1995,ForresterPJ/WitteNS:2006} that we should study
biorthogonal functions with respect to the more general density
\[
\Delta^{(m)}(z;u_0,\dots,u_{2m+5})
=
\frac{\prod\limits_{0\le r<2m+6} \Gampq(u_r z^{\pm 1})}{\Gampq(z^{\pm 2})},
\]
with new balancing condition
\[
\prod_{0\le r<2m+6} u_r = (pq)^{m+1},
\]
and the corresponding contour condition, integrated against the
dif\/ferential
\[
\frac{(p;p)(q;q)}{2}
\frac{dz}{2\pi\sqrt{-1}z}.
\]
Note that if $u_{2m+4}u_{2m+5}=pq$, then the corresponding factors of the
density cancel, and thus we reduce to the order $m-1$ density.  Also, it
will be convenient to multiply the integrands by theta functions, not
elliptic functions; such multiplication has the ef\/fect of shifting the
ba\-lan\-cing condition.  (The extent of the required shift can be determined
via the observation that multiplying a parameter by $q$ multiplies the
integrand by a $p$-theta function; in any event, we will give the explicit
balancing condition for each of the integrals appearing below.)

One natural multivariate analogue of the elliptic beta integral is the
elliptic Selberg integ\-ral~\cite{vanDiejenJF/SpiridonovVP:2000,xforms}, the
higher-order version of which we def\/ine as follows
\begin{gather*}
\II^{(m)}_{n;t;p,q} (u_0,\dots,u_{2m+5}) \\
\qquad{}:=
\frac{(p,p)^n(q;q)^n}{\Gampq(t)^{-n} 2^n n!}
\int_{C^n}
\prod_{1\le i<j\le n}
  \frac{\Gampq(t z_i^{\pm 1}z_j^{\pm 1})}{\Gampq(z_i^{\pm 1}z_j^{\pm 1})}
\prod_{1\le i\le n}
  \frac{\prod\limits_{0\le r<2m+6} \Gampq(u_r z_i^{\pm 1})}{\Gampq(z_i^{\pm 2})}
  \frac{dz_i}{2\pi\sqrt{-1}z_i},
\end{gather*}
where the parameters satisfy the conditions $|t|,|p|,|q|<1$, and
\[
t^{2n-2}\prod_{0\le r<2m+6} u_r = (pq)^{m+1},
\]
and the contour $C$ is chosen so that $C=C^{-1}$, and such that the
interior of $C$ contains every contour of the form $p^i q^j t C$, $i,j\ge
0$, and every point of the form $p^i q^j u_r$, $i,j\ge 0$, $0\le r<2m+6$.
(The latter set of points represents poles of the integrand; if (as often
occurs below) some of these points are not poles, then the corresponding
contour condition can of course be removed.  Similarly, if the cross terms
are holomorphic (e.g., if $t=q$, as is the case below), then $C$ need not
contain the contours $p^i q^j t C$.)  Note that if
$|u_0|$,\dots,$|u_{2m+5}|<1$, then $C$ can be chosen to be the unit circle.
More generally, such a contour exists as long as $p^i q^j t^k u_r u_s$ is
never 1 for $i,j,k\ge 0$, $0\le r,s<2m+6$, and the result is a meromorphic
function on the parameter domain.

When $m=0$, the elliptic Selberg integral can be explicitly evaluated
\cite[Theorem~6.1]{xforms}:
\[
\II^{(0)}_{n;t;p,q}(u_0,u_1,u_2,u_3,u_4,u_5)
=
\prod_{0\le i<n} \Gampq(t^{i+1})\prod_{0\le r<s<6}\Gampq\big(t^i u_r u_s\big),
\]
while the order $1$ elliptic Selberg integral satisf\/ies a transformation
law with respect to the Weyl group $E_7$; more precisely, the renormalized
(holomorphic) function
\begin{equation}
\tilde\II^{(1)}_{n;t;p,q}(u_0,\dots,u_7)
:=
\II^{(1)}_{n;t;p,q}\big(t^{1/2} u_0,\dots,t^{1/2} u_7\big)
\prod_{0\le r<s<8} \Gamma^+_{p,q,t}(t u_r u_s)
\label{eq:eP_tau}
\end{equation}
is invariant under the natural action of $E_7$ on the torus of parameters
\cite[Corollary~9.11]{xforms}.  More importantly for our present purposes, when
$t=q$, the renormalized order $1$ elliptic Selberg integral satisf\/ies an
$E_8$-invariant family of nonlinear dif\/ference equations making it a tau
function for the elliptic Painlev\'e equation \cite[Theorem~5.1]{recur} (for
the relevant def\/inition of tau functions, see~\cite{KajiwaraK/MasudaT/NoumiM/OhtaY/YamadaY:2003}).  As an aside, it
should be noted that \cite{recur} also showed that when $t=q^{1/2}$ or
$t=q^2$, the integral satisf\/ies slightly more complicated analogues of the
tau function identities; as yet, neither a geometric nor an isomonodromy
interpretation of those identities is known.

Since we will be f\/ixing $p$, $q$, and $t=q$ in the sequel, we omit these
parameters from the notation; we will also generally omit $m$, as it can be
determined by counting the arguments.

Consider the following instance of the elliptic Selberg integral:
\[
F_n(x;v) =
x^{-n}v^n\II_n(u_0,\dots,u_{2m+5},qx,pq/x,v,p/v),
\]
satisfying, as usual, the balancing condition
\[
q^{2n-2} \prod_{0\le r<2m+6} u_r = (pq)^{m+1}.
\]
Since
\[
\psi_p(x,z)
:=
x^{-1}\Gampq\big(qx z^{\pm 1},pq/x z^{\pm 1}\big)
=
\frac{\Gampq(qx z^{\pm 1})}{x\Gampq(x z^{\pm 1})}
=
x^{-1}\theta_p\big(x z^{\pm 1}\big),
\]
we see that the integrand of $F_n(x;v)$ is holomorphic in~$x$; indeed, it
dif\/fers from the order~$m$ elliptic Selberg integrand by a factor
\[
\prod_{1\le i\le n} \frac{\psi_p(x,z_i)}{\psi_p(v,z_i)}.
\]
In particular, the $x$-dependent conditions on the contour are irrelevant,
as there are no $x$-dependent poles.  We thus f\/ind that $F_n(x;v)$ is a
$BC_1$-symmetric theta function of degree $n$; that is, it is a holomorphic
function of $x$ satisfying
\begin{gather*}
F_n(1/x;v) =F_n(x;v),\qquad
F_n(px;v)  = \big(px^2\big)^{-n} F_n(x;v).
\end{gather*}
(In general, $BC_n$ denotes the ``hyperoctahedral'' group of signed
permutations, which will act by permutations and taking reciprocals.)  This
function satisf\/ies a form of biorthogonality; to be precise, we have the
following.

\begin{thm}\label{thm:ell_biorth}
Let $G_n(x)$ be any $BC_1$-symmetric theta function of degree $n$,
and let $C$ be any contour satisfying the constraints corresponding to
the parameters $u_0,\dots,u_{2m+5},v,p/v$ with
\[
q^{2n-2} \prod_{0\le r<2m+6} u_r = (pq)^{m+1}.
\]
Then for any $x$ such that the contour $C$ contains $p^i x$ and $p^{i+1}/x$
for all $i\ge 0$,
\begin{gather*}
\frac{(p;p)^2}{2}
\int_C
\frac{F_n(z;v) G_n(z)}{\psi_p(x,z)\psi_p(v,z)}
\Delta(z;u_0,\dots,u_{2m+5})
\frac{dz}{2\pi\sqrt{-1}z}\\
\qquad {}=
G_n(x)
x^{n+1}v^{n+1}
\II_{n+1} (u_0,\dots,u_{2m+5},x,p/x,v,p/v).
\end{gather*}
In particular, if $H_{n-1}(x)$ is a $BC_1$-symmetric theta function of
degree $n-1$, then
\[
\int_C
\frac{F_n(z;v) H_{n-1}(z)}{\psi_p(v,z)}
\Delta(z;u_0,\dots,u_{2m+5})
\frac{dz}{2\pi\sqrt{-1}z}
=
0.
\]
\end{thm}

\begin{proof}
If replace $F_n(z;v)$ by its def\/inition, the result is an $n+1$-dimensional
contour integral over $C^{n+1}$.  Moreover, the integrand is very nearly
symmetric between $z$ and the remaining $n$ integration variables.  To be
precise, if we write the original integration variable as $z_0$, then the
resulting integrand is a $BC_{n+1}$-symmetric factor multiplied by
\[
\frac{G_n(z_0)}{\psi_p(x,z_0)\prod\limits_{1\le i\le n} \psi_p(z_0,z_i)},
\]
which is invariant under the subgroup $BC_1\times BC_n$.  If we average the
integrand over $BC_{n+1}$, this will not change the integral, as the
contour is $BC_{n+1}$-invariant.  We can thus replace the above factor by
the average over cosets:
\[
\frac{1}{n+1}
\sum_{0\le k\le n}
\frac{G_n(z_k)}{\psi_p(x,z_k)\prod\limits_{i\ne k} \psi_p(z_k,z_i)}
=
\frac{1}{n+1}
\frac{G_n(x)}{\prod\limits_{0\le i\le n} \psi_p(x,z_i)};
\]
the identity follows from the fact that if we multiply both sides by
$\prod\limits_{0\le i\le n}\psi_p(x,z_i)$, then both sides are $BC_1$-symmetric
theta functions of degree $n$ in $x$, and agree at the $n+1$ distinct pairs
of points $z_i^{\pm 1}$.

The claim follows immediately.
\end{proof}

\begin{rem}
At the level of orthogonal polynomials, such an $n$-dimensional integral
representation is implicit in \cite{SzegoG:1939}; more precisely, Szeg\"o
gives a representation of orthogonal polynomials as a determinant, but the
Cauchy--Binet identity allows one to turn it into an $n$-dimensional
integral involving the square of a Vandermonde determinant.
\end{rem}

Note that in the above calculation, the $x$-dependent constraint on the
contour was only relevant to the eventual identif\/ication of the
$n+1$-dimensional integral as an elliptic Selberg integral.  We also
observe that if $v$ has the form $u_r/q$, then the parameters $u_r$ and
$p/v$ in the Selberg integrals multiply to $pq$ and thus cancel.  We thus
f\/ind that $F_n(z;u_r/q)$ satisf\/ies biorthogonality with respect to a
general order $m$ instance of $\Delta(z)$.  It will, however, be convenient
to allow general $v$ in the sequel.

We thus see that the integral $F_n(z;v)$ is in some sense an analogue of an
orthogonal polynomial.  Similarly, the $n+1$-dimensional integral of
Theorem \ref{thm:ell_biorth} is analogous to a Cauchy transform of
$F_n(z;v)$, as the integral of $F_n(z;v)$ against a function with a moving
pole.  This suggests that these two integrals should form a row in the
fundamental matrix of our dif\/ference equation.  This leads to the question
of how this row depends on $v$.  Def\/ine
\[
F^+_n(x,v)
=
v^{n+1}x^{n+1}\psi_p(v,x)
\II_{n+1}(u_0,\dots,u_{2m+5},x,p/x,v,p/v),
\]
where the factor $x^{n+1}$ is chosen to make the integrand invariant under
$x\mapsto 1/x$, the factor $v^{n+1}$ for symmetry, and the factor
$\psi_p(v,x)$ to simplify the following identity.

\begin{lem}
The functions $F_n(x;v)$ and $F^+_n(x;v)$ satisfy the identity
\[
F_n(x;v)F^+_n(x,w)
-
F_n(x;w)F^+_n(x,v)
=
\II_n(u_0,\dots,u_{2m+5})F^+_n(v,w).
\]
\end{lem}

\begin{proof}
Taking $G_n(x)=F_n(x;w)$ in Theorem \ref{thm:ell_biorth} gives
\begin{gather*}
F_n(x;w)
F^+_n(x,v)\\
\qquad{}
=
\frac{(p;p)^2}{2}
\int_C
\frac{\psi_p(v,x)}{\psi_p(x,z)\psi_p(v,z)}
F_n(z;v) F_n(z;w)
\Delta(z;u_0,\dots,u_{2m+5})
\frac{dz}{2\pi\sqrt{-1}z}.
\end{gather*}
Thus the two terms on the left-hand side agree except in the f\/irst factors
of the integrands; the dif\/ference of the two integrals can be simplif\/ied
using the addition law, and gives a result independent of $x$; setting
$x=v$ gives the desired result.
\end{proof}

Similarly, we have the following.  Let
\[
F^-_n(x,v)
:=
\psi_p(v,x) x^{1-n}v^{1-n}\II_{n-1}(u_0,\dots,u_{2m+5},qx,pq/x,qv,pq/v).
\]

\begin{lem}
For any $BC_1$-symmetric theta function $G_n$ of degree $n$,
\begin{gather*}
\frac{(p;p)^2}{2}
\int_C
\frac{F^-_n(z,v) G_n(z)\psi_p(x,y)}{\psi_p(x,z)\psi_p(y,z)}
\Delta(z;u_0,\dots,u_{2m+5})
\frac{dz}{2\pi\sqrt{-1}z}\\
\qquad{}=
G_n(x)F_n(v;x)-G_n(y)F_n(v;y).
\end{gather*}
In particular,
\[
F^-_n(x,v)F^+_n(x,w)
+
F_n(v;x)F_n(x;w)
=
\II_n(u_0,\dots,u_{2m+5})F_n(v;w),
\]
and if $H_{n-2}(z)$ is any $BC_1$-symmetric theta function of degree $n-2$,
\[
\int_C
F^-_n(z,v) H_{n-2}(z)
\Delta(z;u_0,\dots,u_{2m+5})
\frac{dz}{2\pi\sqrt{-1}z}
=
0.
\]
\end{lem}

\begin{rem}
In particular, we see that
\[
\frac{F^-_n(z,v)}{\psi_p(v,z)}
\]
is essentially a biorthogonal function of degree $n-1$.
\end{rem}

\begin{thm}
The functions $F^-_n$, $F_n$ and $F^+_n$ satisfy the identities
\begin{gather*}
F^+_n(v,w)
\begin{pmatrix}
F_n(x;u)\\ F^+_n(x,u)
\end{pmatrix}
-
F^+_n(u,w)
\begin{pmatrix}
F_n(x;v)\\
F^+_n(x,v)
\end{pmatrix}
+
F^+_n(u,v)
\begin{pmatrix}
F_n(x;w)\\
F^+_n(x,w)
\end{pmatrix}
 =
0,
 \\
F^+_n(v,w)
\begin{pmatrix}
F^-_n(x,u)\\ -F_n(u;x)
\end{pmatrix}
-
F_n(u;w)
\begin{pmatrix}
F_n(x;v)\\
F^+_n(x,v)
\end{pmatrix}
+F_n(u;v)
\begin{pmatrix}
F_n(x;w)\\
F^+_n(x,w)
\end{pmatrix}
 =0,
 \\
F_n(v;w)
\begin{pmatrix}
F^-_n(x,u)\\ -F_n(u;x)
\end{pmatrix}
-
F_n(u;w)
\begin{pmatrix}
F^-_n(x,v)\\ -F_n(v;x)
\end{pmatrix}
+F^-_n(u,v)
\begin{pmatrix}
F_n(x;w)\\
F^+_n(x,w)
\end{pmatrix}
 =0,
 \\
F^-_n(v,w)
\begin{pmatrix}
F^-_n(x,u)\\ -F_n(u;x)
\end{pmatrix}
-
F^-_n(u,w)
\begin{pmatrix}
F^-_n(x,v)\\ -F_n(v;x)
\end{pmatrix}
+
F^-_n(u,v)
\begin{pmatrix}
F^-_n(x,w)\\ -F_n(w;x)
\end{pmatrix}
 =0.
\end{gather*}
\end{thm}

\begin{proof}
Each identity is the Pl\"ucker relation for the $2\times 3$ matrix formed
by concatenating the three column vectors that appear.  In the f\/irst two
cases, we have already computed the requisite minors; the remaining minor
follows as a special case of the third identity, which can be derived by
eliminating a common term from two instances of the second identity.
\end{proof}

\begin{rem}
Note that the proof of these identities didn't require the balancing
condition, or even that the biorthogonality density was $\Delta$.
Furthermore, the only way in which the proof depended on properties of
elliptic functions was in the fact that $\psi_p$ satisf\/ies a partial
fraction decomposition result.  If we generalize the results with this in
mind, we f\/ind that these are precisely the generalized Fay
identities of \cite{AdlerM/vanMoerbekeP:1999,recur}.

We also note that the change of basis from $F_n$ to $F^-_n$ can be
interpreted as relating degree $n$ biorthogonal functions to degree $n-1$
biorthogonal functions; i.e., the analogue of the three-term recurrence for
orthogonal polynomials.
\end{rem}

We thus see that, as functions of $x$, the vectors
\[
\begin{pmatrix}
F_n(x;v)\\
F^+_n(x,v)
\end{pmatrix}
\]
and
\[
\begin{pmatrix}
F^-_n(x,v)\\ -F_n(v;x)
\end{pmatrix}
\]
for all $v\in \C^*$, together span only a $2$-dimensional space, and the
change of basis matrix between any two such bases of this $2$-dimensional
space is computable in terms of $F_n$, $F^{\pm}_n$.  And, naturally, the
choice of basis will have no ef\/fect on the resulting dif\/ference equation
beyond conjugation by a matrix independent of $x$.  Since it will be useful
to allow fairly general bases, we extend the notation by def\/ining values
for $F^+_n$ on hatted arguments (equivalently, def\/ining $F^+_n$ as a
function on $(\C^*\disj \C^*)^2$, where $\C^*\disj \C^*$ is the disjoint
union of two copies of $\C^*$; thus a given element $v\in \C^*$ corresponds
to two elements of $\C^*\disj \C^*$, denoted by $v$ and $\hat{v}$
respectively), as follows:
\begin{alignat*}{3}
& F^+_n(\hat{v},w)  := F_n(v;w),\qquad &&
F^+_n(v,\hat{w})  := -F_n(w;v),& \\
& F^+_n(\hat{v},\hat{w})  :=F^-_n(v,w),\qquad &&
F_n(x;\hat{v}) := F^-_n(x,v);&
\end{alignat*}
note that this extension of $F^+_n$ preserves its antisymmetry.  Note
that in this notation, the identities relating $F_n$, $F^\pm_n$ reduce to
the single identity
\[
F^+_n(w,x)F^+_n(y,z)
-F^+_n(w,y)F^+_n(x,z)
+F^+_n(w,z)F^+_n(x,y)
=0,
\]
for all $w,x,y,z\in \C^*\disj \C^*$, and the minors used in the Pl\"ucker
relations follow from the special case
\[
F^+_n(\hat{v},v) = F_n(v;v) = \II_n(u_0,\dots,u_{2m+5}).
\]

To proceed further, we will need to understand how our functions behave
under the monodromy action $x\mapsto px$; it will also turn out to be
useful to know how $x\mapsto 1/x$ acts.  Easiest of all is $x\mapsto p/x$;
in that case, the elliptic Selberg integral itself is manifestly invariant,
so we need simply consider how the prefactors transform:
\begin{gather*}
F_n(p/x;v)    = \big(x^2/p\big)^n F_n(x;v),\qquad
F^+_n(p/x,v)  = \big(p/x^2\big)^n F^+_n(x,v).
\end{gather*}
For $x\mapsto 1/x$, we similarly have
\[
F_n(1/x;v) = F_n(x;v).
\]
However, for $F^+_n(1/x,v)$, while the integrand remains constant, the
constraints on the contour change.  Assume for the moment that $v\in \C^*$,
and choose a $BC_1$-symmetric theta function $G_n(z)$ of degree $n$ such
that $G_n(x)\ne 0$, so that
\[
F^+_n(x,v)
=
\frac{(p;p)^2}{2}
\int_C
\frac{\psi_p(v,x) F_n(z;v) G_n(z)}{\psi_p(x,z)\psi_p(v,z)G_n(x)}
\Delta(z;u_0,\dots,u_{2m+5})
\frac{dz}{2\pi\sqrt{-1}z}.
\]
Then $x\mapsto 1/x$ leaves the integrand the same, but moves the contour
through $x$ and $1/x$.  Thus $F^+_n(1/x,v)-F^+_n(x,v)$ can be computed by
residue calculus; by symmetry, we f\/ind that it is twice the residue at
$z=1/x$:
\[
F^+_n(1/x,v)-F^+_n(x,v)
=
x^{-1}\theta_q\big(x^2\big) F_n(x;v)
\prod_{0\le r<2m+6} \Gampq\big(u_r x^{\pm 1}\big).
\]
Putting this together, we obtain the following.

\begin{lem}
The functions $F_n$ and $F^+_n$ have the monodromy action
\[
\rowvec{F_n(1/x;v)}{F^+_n(1/x,v)}
\begin{matrix}{}={}\\{}\end{matrix}
\rowvec{F_n(x;v)}{F^+_n(x,v)}
\begin{pmatrix}
1 & x^{-1}\theta_q(x^2)
\prod\limits_{0\le r<2m+6} \Gampq(u_r x^{\pm 1})\vspace{1mm}\\
0 & 1
\end{pmatrix}
\]
and
\begin{gather*}
 \begin{pmatrix}F_n(px;v) & F^+_n(px,v)\end{pmatrix} \\
\qquad\begin{matrix}{}={}\\{}\end{matrix}
\rowvec{F_n(x;v)}{F^+_n(x,v)}
\begin{pmatrix}
(px^2)^{-n} & (px^2)^n x^{-1}\theta_q(x^2)
\prod\limits_{0\le r<2m+6} \Gampq(u_r x^{\pm 1})\vspace{1mm}\\
0 & (px^2)^n
\end{pmatrix}
\end{gather*}
valid for all $v\in \C^*\disj \C^*$.
\end{lem}

\begin{proof}
The only thing to check is that it extends to the other copy of $\C^*$, but
this follows immediately from the facts that the monodromy is independent of
$v\in \C^*$, and that for all $v\in \C^*\disj \C^*$, the row vectors lie in
the same $2$-dimensional space.
\end{proof}

This is not quite a $q$-theta $p$-dif\/ference equation as we would wish, but
it is straightforward to turn it into a $q$-theta $p$-dif\/ference equation.
Def\/ine a $2\times 2$ meromorphic matrix $M_n(z;v,w)$ for $v,w\in \C^*\disj
\C^*$:
\[
M_n(z;v,w)
:=
\begin{pmatrix}
F_n(z;v) & z^{-1}\theta_p(z^2)F^+_n(z,v)/\Delta(z;u_0,\dots,u_{2m+5})\vspace{1mm}\\
F_n(z;w) & z^{-1}\theta_p(z^2)F^+_n(z,w)/\Delta(z;u_0,\dots,u_{2m+5})
\end{pmatrix}.
\]

\begin{thm}
The matrix $M_n(z;v,w)$ is a meromorphic fundamental matrix for a $p$-theta
$q$-difference equation with multiplier $q^{-2n}$.  The isomonodromy
class of the equation is independent of $v$ and $w$, and invariant under
all permutations of the parameters and all shifts
\[
(u_0,\dots,u_{2m+5})\mapsto \big(q^{k_0}u_0,\dots,q^{k_{2m+5}}u_{2m+5}\big)
\]
with $k_r\in \Z$ such that
\[
\sum_{0\le r<2m+6} k_r = 0;
\]
and invariant under simultaneous negation of all parameters.
The theta-isomonodromy class is further invariant under all shifts
\[
(u_0,\dots,u_{2m+5},z,n)\mapsto
\big(q^{k_0}u_0,\dots,q^{k_{2m+5}}u_{2m+5},q^l z,n+\nu\big)
\]
with $l\in \frac{1}{2}\Z$, $k_r\in l+\Z$, $\nu\in \Z$ such that
\[
2\nu + \sum_{0\le r<2m+6} k_r = 0.
\]
In addition, the associated shift matrix $A(z)$ satisfies the symmetry
$A(1/qz)A(z)=1$.
\end{thm}

\begin{proof}
Most of the claims follow immediately from the fact that $M(z;v,w)$
satisf\/ies the $q$-theta $p$-dif\/ference equation
\begin{gather*}
M_n(pz;v,w)\\
\qquad{} =
M_n(z;v,w)
\begin{pmatrix}
(pz^2)^{-n}&
(pz^2)^{n-2} \Delta(z;u_0,\dots,u_{2m+5})/\Delta(pz;u_0,\dots,u_{2m+5})
\vspace{1mm}\\
0 &-(pz^2)^{n-2} \Delta(z;u_0,\dots,u_{2m+5})/\Delta(pz;u_0,\dots,u_{2m+5})
\end{pmatrix},
\end{gather*}
and has determinant
\[
\det(M_n(z;v,w))
=
\II_n(u_0,\dots,u_{2m+5})F^+_n(v,w)
\frac{z^{-1}\theta_p(z^2)}{\Delta(z;u_0,\dots,u_{2m+5})},
\]
so is nonsingular.  The only additional thing to check for the isomonodromy
claims is that
\[
\big(pz^2\big)^{2n-2} \Delta(z;u_0,\dots,u_{2m+5})/\Delta(pz;u_0,\dots,u_{2m+5})
\]
is invariant 
under all of the
stated transformations.  The symmetry of $A(z)$ follows immediately from
the symmetry
\begin{equation*}
M(1/z;v,w)
=
M(z;v,w)
\begin{pmatrix}
1 & 1 \\
0 & -1
\end{pmatrix}.
\tag*{\qed}
\end{equation*}
\renewcommand{\qed}{}
\end{proof}

\begin{rem}
One can avoid the appearance of theta-isomonodromy above at the cost of
introducing some ``apparent'' singularities, and an additional parameter
controlling the location of those singularities.  Indeed, if one def\/ines
\[
M'_n(z;x;v,w)
:=
\frac{\Gampq(x z^{\pm 1})}
     {\Gampq(q^n x z^{\pm 1})}
M_n(z;v,w),
\]
then $M'_n$ satisf\/ies the same transformation law as $M$ with respect to
$z\mapsto 1/z$, while under $z\mapsto pz$, one has
\begin{gather*}
M'_n(pz;x;v,w)\\
\qquad{}
=
M'_n(z;x;v,w)
\begin{pmatrix}
1&(pz^2)^{2n-2} \Delta(z;u_0,\dots,u_{2m+5})/\Delta(pz;u_0,\dots,u_{2m+5})
\vspace{1mm}\\
0 &-(pz^2)^{2n-2} \Delta(z;u_0,\dots,u_{2m+5})/\Delta(pz;u_0,\dots,u_{2m+5})
\end{pmatrix}.
\end{gather*}
Thus the associated shift matrix
\[
A'_n(z;x;v,w)
=
\frac{\theta_p(xz,q^{n-1}x/z)}
     {\theta_p(q^n xz,x/qz)}
A_n(z;v,w)
\]
is {\em elliptic}, with the same symmetry as $A$, and every shift
\[
(u_0,\dots,u_{2m+5},z,n,x)\mapsto
\big(q^{k_0}u_0,\dots,q^{k_{2m+5}}u_{2m+5},q^l z,n+\nu,q^{l'}x\big)
\]
with $l\in \frac{1}{2}\Z$, $l',k_r\in l+\Z$, $\nu\in \Z$ such that
\[
2\nu + \sum_{0\le r<2m+6} k_r = 0
\]
gives rise to a true isomonodromy transformation of this elliptic
dif\/ference equation, with associated operator
\[
B'(z;x;v,w)
=
\frac{\Gampq(q^{l'+l}x z,q^{l'-l}x/z)}
     {\Gampq(x z,x/z)}
\frac{\Gampq(q^{n} x z,q^{n} x/z)}
     {\Gampq(q^{n+\nu+l+l'} x z,q^{n+\nu+l'-l} x/z)}
B(z;v,w).
\]
In particular, the isomonodromy transformations dif\/fer from the
corresponding theta-isomo\-nodromy transformations by a meromorphic
theta function factor depending only on $\nu$, $l'$, $l$.  There is also
an isomonodromy transformation between $A'_n(z;x;v,w)$ and
$A'_n(z;x';v,w)$, for arbitrary~$x'$, but the corresponding $B$ matrix is
(generically) multiplication by an elliptic function of degree~$2n$.  It
follows that only those parameter shifts satisfying the integrality
condition above can extend to arbitrary solutions of the elliptic
Painlev\'e equation (for which one ef\/fectively has noninteger~$n$).
\end{rem}

\begin{rem}
The lattice of possible $k_r$ vectors is the lattice $D_{2m+6}^+$ obtained
by adjoining the vector
\[
(1/2,1/2,\dots,1/2,-1/2,-1/2,\dots,-1/2)
\]
of sum 0 to the root lattice $D_{2m+6}$.  In particular, when $m=1$,
this lattice is precisely the root lattice $E_8$.
\end{rem}

\begin{rem}
  The symmetry of $A(z)$ is precisely the condition for the pair of
  equations
\[
v(qz)=A(z)v(z),\qquad  v(1/z)=v(z)
\]
to be formally consistent.  It then follows from a dif\/ferent special case
of \cite[Theorem~3]{PraagmanC:1986} that there exists a nonsingular
meromorphic matrix $\hat M$ such that
\[
\hat{M}(qz)=A(z)\hat{M}(z),\qquad  \hat{M}(1/z)=\hat{M}(z).
\]
We can give such a matrix explicitly, for instance
\[
\hat{M}(z) = M(z) \begin{pmatrix} 1&-1/2\\0&1\end{pmatrix}
\begin{pmatrix} 1&0\vspace{1mm}\\0& \dfrac{\theta_p(az^{\pm 1},bz^{\pm
      1})}{z^{-1}\theta_p(z^2)}
\end{pmatrix}.
\]
\end{rem}

\section{The dif\/ference equation}\label{section4}

Naturally, simply knowing the {\em existence} of a dif\/ference equation with
associated isomonodromy transformations is of strictly limited usefulness,
so we would like to be more explicit about the equation, and at the very
least generators of the group of monodromy-preserving transformations.

The f\/irst thing we will need to understand about the shift matrix is the
locations of its singularities; i.e., the points where the coef\/f\/icients
have poles or the determinant has a zero.  This in turn depends on
determining the polar divisor of $F^+_n$.  Def\/ine
\[
(x;p,q):=\prod_{0\le i,j} (1-p^i q^j x),
\]
with the usual multiple argument conventions.

\begin{lem}
The function
\[
\left(\prod_{0\le r<2m+6} (u_r x,pu_r/x;p,q)\right)
F^+_n(x,v)
\]
is holomorphic for $x\in \C^*$.  If $v\in \C^*$, the function vanishes at
$x=v$, $x=p/v$.
\end{lem}

\begin{proof}
As before, assuming $v\in \C^*$, we have
\begin{gather*}
F^+_n(x,v)
 =
\frac{(p;p)^2}{2}
\int_C
\frac{\psi_p(v,x) F_n(z;v) G_n(z)}{\psi_p(x,z)\psi_p(v,z) G_n(x)}
\Delta(z;u_0,\dots,u_{2m+5})
\frac{dz}{2\pi\sqrt{-1}z} \\
\phantom{F^+_n(x,v)}{}
=
\frac{(p;p)^2}{2}
\int_C
\frac{xv\psi_p(v,x) F_n(z;v) G_n(z)}{G_n(x)}
\Delta(z;u_0,\dots,u_{2m+5},x,p/x,v,p/v)
\frac{dz}{2\pi\sqrt{-1}z}
.
\end{gather*}
where $G_n(z)$ is any $BC_1$-symmetric $p$-theta function of degree $n$ not
vanishing at~$x$.  Since $F_n(z;v)$ is holomorphic in $z$, we may apply
Lemma~10.4 of \cite{xforms} (note that condition 3 of that lemma reduces in
our case to the balancing condition) to conclude that
\begin{gather*}
\left(\prod_{0\le r<2m+6}(u_r x,u_r p/x;p,q)\right)
\big(xv,px/v,pv/x,p^2/xv;p,q\big)
(xv)^{-1}\psi_p(v,x)^{-1}
G_n(x)
F^+_n(x,v)
\end{gather*}
is holomorphic in $x$.  (The conclusion concerning the $x$-independent
poles is not useful to us, as~$F_n(z;v)$ certainly has singularities that
depend on the remaining parameters.)  This nearly gives us the desired
result, except for the factor~$G_n(x)$, which disappears by the fact that
$F^+_n(x,v)$ is independent of~$G_n$, and the additional factor
\[
\big(qxv,pqx/v,pqv/x,p^2q/xv;p,q\big).
\]
This latter factor can be eliminated, and the result extended to $v\in
\C^*\disj \C^*$, by expressing $F^+_n(z,v)$ as a linear combination of
$F^+_n(z,w)$ and $F^+_n(z,w')$, which for generic $w$ and $w'\in \C^*$ are
holomorphic at the of\/fending points.
\end{proof}

\begin{thm}
The matrix
\begin{gather*}
\tilde{A}_n(z;v,w)
 :=
\left(q^{-1} z^{-2} \prod_{0\le r<2m+6} \theta_p(u_r z)\right)
A_n(z;v,w) \\
\phantom{\tilde{A}_n(z;v,w)}{} =
\left(q^{-1} z^{-2} \prod_{0\le r<2m+6} \theta_p(u_r z)\right)
M_n(qz;v,w)M_n(z;v,w)^{-1}
\end{gather*}
is holomorphic in $z$, with determinant
\[
\det\big(\tilde{A}_n(z;v,w)\big)
=
\prod_{0\le r<2m+6} \theta_p(u_r z,u_r/qz),
\]
satisfies the $p$-theta transformation law
\[
\tilde{A}_n(pz;v,w)
=
\big(pqz^2\big)^{-m-3}
\tilde{A}_n(z;v,w)
,
\]
and has the symmetry
\[
\tilde{A}_n(1/qz;v,w)
=
\begin{pmatrix}
 0 & -1\\
 1 &  0
\end{pmatrix}
\tilde{A}_n(z;v,w)^t
\begin{pmatrix}
 0 & 1\\
-1 & 0
\end{pmatrix}.
\]
\end{thm}

\begin{proof}
  The formula for $\det(\tilde{A}_n(z;v,w))$ follows immediately from the
  formula for the determinant $\det(M_n(z;v,w))$.  Similarly, the fact that
  $A_n(1/qz;v,w)A_n(z;v,w)=1$ becomes
\[
\tilde{A}_n(1/qz;v,w)\tilde{A}_n(z;v,w) = \det\tilde{A}_n(z;v,w);
\]
the symmetry of $\tilde{A}_n(z;v,w)$ follows from the usual formula for the
inverse of a $2\times 2$ matrix:
\[
C^{-1} =
\det(C)^{-1}
\begin{pmatrix}
C_{22} & -C_{12}\\
-C_{21} & C_{11}
\end{pmatrix}.
\]

Another use of this formula allows us to explicitly write down the inverse of
$M_n(z;v,w)$.  This, in turn, allows us to express the entries of
$\tilde{A}(z)$ as polynomials in $F_n$ and $F^+_n$ with coef\/f\/icients that
are (holomorphic) $p$-theta functions in $z$:
\begin{gather*}
\tilde{A}_n(z;v,w)
 =
\begin{pmatrix}
F_n(qz;v) & F^+_n(qz,v)\\
F_n(qz;w) & F^+_n(qz,w)
\end{pmatrix}
\begin{pmatrix} a(z) &0\\0&b(z)\end{pmatrix}
\begin{pmatrix}
F^+_n(z,w) & -F^+_n(z,v)\\
-F_n(z;w) & F_n(z;v)
\end{pmatrix} \\
\phantom{\tilde{A}_n(z;v,w)}{} =
\begin{pmatrix}
F_n(qz;v) & b(z)F^+_n(qz,v)\\
F_n(qz;w) & b(z)F^+_n(qz,w)
\end{pmatrix}
\begin{pmatrix}
a(z) F^+_n(z,w) & -a(z)F^+_n(z,v)\\
-F_n(z;w) & F_n(z;v)
\end{pmatrix},
\end{gather*}
where
\begin{gather*}
a(z)  =
\frac{q^{-1} z^{-2} \prod\limits_{0\le r<2m+6} \theta_p(u_r z)}
     {\II_n(u_0,\dots,u_{2m+5})F^+_n(v,w)},
 \qquad
b(z) =
\frac{q z^2 \prod\limits_{0\le r<2m+6} \theta_p(u_r/qz)}
     {\II_n(u_0,\dots,u_{2m+5})F^+_n(v,w)}.
\end{gather*}

In particular, the only possible poles of $\tilde{A}_n$ come from poles of
$a(z)F^+_n(z,v)$, $a(z)F^+_n(z,w)$, $b(z)F^+_n(qz,v)$ and
$b(z)F^+_n(qz,w)$, or, equivalently, poles of
\[
a(z)\prod_{0\le r<2m+6} (u_r z,pu_r/z;p,q)^{-1}
\sim
\prod_{0\le r<2m+6}
\frac{(p/u_rz;p)}
{(u_r qz,pu_r/z;p,q)}
\]
(where $\sim$ here and below denotes that the two functions have the same
zeros and poles) and
\[
b(z)\prod_{0\le r<2m+6} (u_r qz,pu_r/qz;p,q)^{-1}
\sim
\prod_{0\le r<2m+6}
\frac{(qz/u_r;p)}
{(u_r qz,pu_r/z;p,q)}.
\]
It follows that
\[
\prod_{0\le r<2m+6} (u_r qz,pu_r/z;p,q)\tilde{A}_n(z;v,w)
\]
is holomorphic in $z$.  But the entries of $\tilde{A}_n(z;v,w)$ are
meromorphic $p$-theta functions, and thus their divisors are periodic in
$p$.  Since the remaining set of potential poles contains no $p$-periodic
subset, there are in fact no surviving poles.
\end{proof}

We can also compute the value of $\tilde{A}_n(z;v,w)$ at a number of
points.

\begin{thm}
The matrix $\tilde{A}_n(z;v,w)$ has the special values
\begin{gather*}
 \tilde{A}_n(u_s/q;v,w)
 =
\frac{qu_s^{-2} \prod\limits_{0\le r<2m+6} \theta_p(u_ru_s/q)}
     {\II_n(u_0,\dots,u_{2m+5})F^+_n(v,w)}
\begin{pmatrix}
F_n(u_s;v) \\
F_n(u_s;w)
\end{pmatrix}
\rowvec{F^+_n(u_s/q,w)}{-F^+_n(u_s/q,v)}
, \\
 \tilde{A}_n(1/u_s;v,w)
 =
\frac{qu_s^{-2} \prod\limits_{0\le r<2m+6} \theta_p(u_ru_s/q)}
     {\II_n(u_0,\dots,u_{2m+5})F^+_n(v,w)}
\begin{pmatrix}
F^+_n(u_s/q,v)\\
F^+_n(u_s/q,w)
\end{pmatrix}
\rowvec{{-}F_n(u_s;w)}{F_n(u_s;v)}
,
\end{gather*}
for $0\le s<2m+6$.  In addition, we have the four values
\begin{gather*}
\tilde{A}_n\big(q^{-1/2};v,w\big)
=
\prod_{0\le r<2m+6} \theta_p\big(u_r q^{-1/2}\big)
\begin{pmatrix} 1&0\\0&1\end{pmatrix}
,
\\
\tilde{A}_n\big({-}q^{-1/2};v,w\big)
=
\prod_{0\le r<2m+6} \theta_p\big({-}u_r q^{-1/2}\big)
\begin{pmatrix} 1&0\\0&1\end{pmatrix}
,
\\
\tilde{A}_n\big((p/q)^{1/2};v,w\big)
=
q^{-n} p^{-1} \prod_{0\le r<2m+6} \theta_p\big(u_r (p/q)^{1/2}\big)
\begin{pmatrix} 1&0\\0&1\end{pmatrix}
,
\\
\tilde{A}_n\big({-}(p/q)^{1/2};v,w\big)
=
q^{-n} p^{-1} \prod_{0\le r<2m+6} \theta_p\big({-}u_r (p/q)^{1/2}\big)
\begin{pmatrix} 1&0\\0&1\end{pmatrix}
\end{gather*}
at the ramification points $($fixed points of $z\mapsto 1/qz$ modulo $\langle
p\rangle)$.
\end{thm}

\begin{proof}
We f\/irst observe that at $z=u_s/q$, $b(z)F^+_n(qz,v)$ and $b(z)F^+_n(qz,w)$
vanish, and thus the formula for $\tilde{A}(u_s/q)$ simplif\/ies as stated.
The second set of special values follows similarly from the vanishing of
$a(z)F^+_n(z,v)$ and $a(z)F^+_n(z,w)$ at $z=p/u_s$, together with the
$p$-theta law of $\tilde{A}$.

When $z=\pm q^{-1/2}$, so that $qz=1/z$, we f\/ind $a(z)=b(z)$,
$F_n(qz,v)=F_n(z,v)$, and $F^+_n(qz,v)=F^+_n(z,v)$; the last dif\/ference
vanishes due to the factor $\theta_q(z^2)=0$ in the relevant residue.  The
expression for $\tilde{A}(\pm q^{-1/2})$ thus simplif\/ies immediately.
Similarly, at $z=\pm \sqrt{p/q}$, we have $qz=p/z$, and again the entries
immediately simplify.
\end{proof}

Note that the symmetry of $\tilde{A}$ and the elementary values at
the ramif\/ication points imply that the matrix is already determined
by its values at $u_s/q$ for any $m+2$ values of $s$ (assuming $u_s$
are generic); the above special values are thus highly overdetermined.
It is also worth noting that if $v$ and $w$ are of the form $u_r/q$
or $\widehat{u_r}$ (with dif\/ferent values of $r$), then all but two
pairs of special values can be expressed entirely in terms of order
$m$ elliptic Selberg integrals with shifted parameters.  This
specialization also has the ef\/fect of causing the kernel of
$\tilde{A}_n(u_r/q)$ and image of $\tilde{A}_n(1/u_r)$ (or vice
versa, as appropriate) to be coordinate vectors.  In particular,
it follows that cross-ratios of kernel vectors are themselves ratios
of order $m$ elliptic Selberg integrals, so long as the four arguments
in question do not contain both $u_s/q$ and $1/u_s$ for some $s$.  For
instance, the cross-ratio of the vectors $\ker\tilde{A}_n(u_0/q),\dots,\ker\tilde{A}_n(u_3/q)$ is
\[
\frac{F^+_n(u_0/q,u_2/q)F^+_n(u_1/q,u_3/q)}{F^+_n(u_0/q,u_3/q)F^+_n(u_1/q,u_2/q)}.
\]

We also note that when $n=0$, $F_0(z;\hat{w})=F^-_0(z,w)=0$, and thus if
$v$ is not ``hatted'', then $\tilde{A}_0(z;v,\hat{w})$ is well-def\/ined and
triangular; note in particular that $F^+_0(v,\hat{w})=1$.  In particular,
it follows that for any~$n\ge 0$, $\tilde{A}_n(z;v,w)$ is isomonodromic to
a triangular shift operator with at most the same number of singularities.

It remains to consider the isomonodromy transformations.  Changing $v$ and
$w$ is straightforward, as we have already observed; the precise
isomonodromy transformation follows from
\begin{equation}
M_n(z;v',w')
=
\frac{1}{F^+_n(v,w)}
\begin{pmatrix}
F^+_n(v',w) & -F^+_n(v',v)\\
F^+_n(w',w) & -F^+_n(w',v)
\end{pmatrix}
M_n(z;v,w).
\label{eq:isomono_vw_vpwp}
\end{equation}
With this in mind, we can feel free to make choices for $v$ and $w$ if this
will simplify the expressions for the remaining isomonodromy
transformations.

We f\/irst consider the case of integer shifts, or in other words shifts
under the lattice $D_{2m+6}$.  It thus suf\/f\/ices to consider the two
cases $(u_0,u_1,n)\mapsto (qu_0,u_1/q,n),(qu_0,qu_1,n-1)$.

\begin{lem}
We have the $($isomonodromy$)$ transformations
\begin{gather}
M_n(z;u_0,\widehat{u_1/q};qu_0,u_1/q,u_2,\dots,u_{2m+5})\nonumber\\
\qquad{} =
(u_0q/u_1)^n
\begin{pmatrix}
1 & 0\\
0 & \dfrac{\theta_p(u_1 z^{\pm 1}/q)}{\theta_p(u_0 z^{\pm 1})}
\end{pmatrix}
M_n(z;u_1/q,\widehat{u_0};u_0,\dots,u_{2m+5})\label{eq:isomono_ud}
\end{gather}
and
\begin{gather}
M_{n-1}(z;u_0,u_1;qu_0,qu_1,u_2,\dots,u_{2m+5})\nonumber\\
\qquad{} =
(u_0u_1)^{n-1}
\begin{pmatrix}
\dfrac{u_1}{\theta_p(u_1 z^{\pm 1})}&0\\
0&\dfrac{u_0}{\theta_p(u_0 z^{\pm 1})}
\end{pmatrix}
M_n(z;\widehat{u_1},\widehat{u_0};u_0,\dots,u_{2m+5}).\label{eq:isomono_uu}
\end{gather}
\end{lem}

\begin{rem}\sloppy
In the version with apparent singularities controlled by $x$, the
transforma\-tions~\eqref{eq:isomono_vw_vpwp} and~\eqref{eq:isomono_ud} remain
unchanged, while~\eqref{eq:isomono_uu} should be multiplied by~$\theta_p(q^{n-1}x z^{\pm 1})$.
\end{rem}

It remains to consider the case of a half-integer shift, say
\begin{gather*}
(u_0,\dots,u_{m+2},u_{m+3},\dots,u_{2m+5},z,v,w) \\
\qquad{} \mapsto
\big(q^{1/2}u_0,\dots,q^{1/2}u_{m+2},
 q^{-1/2}u_{m+3},\dots,q^{-1/2}u_{2m+5},q^{1/2}z,q^{1/2}v',q^{1/2}w'\big).
\end{gather*}
In this case, just as with $A$ itself, we cannot give a closed-form
expression, so must be content with a description of the singularities and
a suf\/f\/iciently large set of special values.  Let $B_n(z;v,w;v',w')$ denote
this isomonodromy transformation; that is,
\[
B_n(z;v,w;v',w')
=
M_n\big(q^{1/2}z;q^{1/2}v',q^{1/2}w';u'_0,\dots,u'_{2m+5}\big)
M_n(z;v,w;u_0,\dots,u_{2m+5})^{-1},
\]
where
\[
u'_r = \begin{cases} q^{1/2}u_r, & 0\le r<m+3,\\
                     q^{-1/2}u_r, & m+3\le r<2m+6.
\end{cases}
\]
Note that
\[
B_n(1/qz;v,w;v',w')^{-1} B_n(z;v,w;v',w')
=
A_n(z;v,w).
\]

For convenience, we def\/ine
\[
G^+_n(v,w):=F^+_n\big(q^{1/2}v,q^{1/2}w;u'_0,\dots,u'_{2m+5}\big),
\]
and similarly for $G_n(v;w)$; this will free us to again omit most of the
parameters from the functions $F_n$ and $F^+_n$.

\begin{thm}
The matrix
\[
\tilde{B}_n(z;v,w;v',w')
:=
\left(
\big(q^{1/2}z\big)^{-1}
\prod_{0\le r<m+3} \theta_p(u_r z)\right)
B_n(z;v,w;v',w')
\]
is holomorphic in $z$ with $p$-theta law
\[
\tilde{B}_n(pz;v,w;v',w')
=
\left(
q^{-n} p^{-1}
(-1)^{m+3}
\prod_{0\le r<m+3} u_r^{-1}
\right)
z^{-m-3}
\tilde{B}_n(z;v,w;v',w'),
\]
has determinant
\[
-q^{-1/2}
\frac{\II_n(u'_0,\dots,u'_{2m+5}) G^+_n(v',w')}
{\II_n(u_0,\dots,u_{2m+5}) F^+_n(v,w)}
\prod_{0\le r<m+3} \theta_p(u_r z)
\prod_{m+3\le r<2m+6} \theta_p(u_r/qz),
\]
and has the special values
\begin{gather*}
\tilde{B}_n(1/u_s;v,w;v',w')\\
\qquad{}
=
\frac{
-u_s^{-1}
\prod\limits_{m+3\le r<2m+6}\theta_p(u_ru_s/q)}
{\II_n(u_0,\dots,u_{2m+5})
F^+_n(v,w)
}
\begin{pmatrix}
G^+_n(u_s/q,v')\\
G^+_n(u_s/q,w')
\end{pmatrix}
\rowvec{{-}F_n(u_s;w)}{F_n(u_s;v)}
,
\end{gather*}
for $0\le s<m+3$, and
\begin{gather*}
\tilde{B}_n(u_s/q;v,w;v',w') \\
\qquad{}=
\frac{(q^{-1/2}u_s)^{-1}\prod\limits_{0\le r<m+3} \theta_p(u_r u_s/q)}
{\II_n(u_0,\dots,u_{2m+5}) F^+_n(v,w)}
\begin{pmatrix}
G_n(u_s/q;v')\\
G_n(u_s/q;w')
\end{pmatrix}
\rowvec{F^+_n(u_s/q,w)}{-F^+_n(u_s/q,v)}
\end{gather*}
for $m+3\le s<2m+6$.
\end{thm}

\begin{proof}
As in the computation for $\tilde{A}_n$, we can use the known determinant
of $M_n$ to write the entries of $\tilde{B}_n$ in terms of $F_n$, $F^+_n$,
$G_n$, and $G^+_n$. We f\/ind
\[
\tilde{B}_n(z;v,w;v',w')
=
\begin{pmatrix}
G_n(z;v') & d(z)G^+_n(z,v')\\
G_n(z;w') & d(z)G^+_n(z,w')
\end{pmatrix}
\begin{pmatrix}
c(z)F^+_n(z,w) & -c(z)F^+_n(z,v)\\
-F_n(z;w) & F_n(z;v)
\end{pmatrix},
\]
where
\begin{gather*}
c(z)  = \frac{(q^{1/2}z)^{-1}\prod\limits_{0\le r<m+3} \theta_p(u_r z)}{\II_n(u_0,\dots,u_{2m+5}) F^+_n(v,w)},\qquad
d(z)  = \frac{-z\prod\limits_{m+3\le r<2m+6} \theta_p(u_r/qz)}{\II_n(u_0,\dots,u_{2m+5}) F^+_n(v,w)}.
\end{gather*}
Since
\begin{gather*}
c(z)F^+_n(z,v)
 \sim
\frac{\prod\limits_{0\le r<m+3} (p/u_rz;p)}
     {\prod\limits_{0\le r<m+3} (q u_r z,pu_r/z;p,q)\prod\limits_{m+3\le r<2m+6} (u_r z,pu_r/z;p,q)},
\\
d(z)G^+_n(z,v')
 \sim
\frac{\prod\limits_{m+3\le r<2m+6} (qz/u_r;p)}
     {\prod\limits_{0\le r<m+3} (q u_r z,pu_r/z;p,q)\prod\limits_{m+3\le r<2m+6} (u_r z,pu_r/z;p,q)},
\notag
\end{gather*}
and the coef\/f\/icients of $\tilde{B}_n$ are $p$-theta functions, we conclude
as before that $\tilde{B}_n$ is holomorphic.

The special values again follow by choosing $z$ so that $c(z)$ or $d(z)$
vanishes, and using the $p$-theta law as appropriate.
\end{proof}

\begin{rem}
The relation between $A$ and $B$ becomes, via the usual expression for the
inverse, the expression
\[
\begin{pmatrix}0 & -1\\1 & 0\end{pmatrix}
\tilde{B}_n(1/qz;v,w;v',w')^t
\begin{pmatrix}0 & 1\\-1 & 0\end{pmatrix}
\tilde{B}_n(z;v,w;v',w')
=
C
\tilde{A}_n(z;v,w),
\]
where
\[
C=
\frac{
-q^{-1/2}
\II_n(u'_0,\dots,u'_{2m+5})G^+_n(v',w')
}{
\II_n(u_0,\dots,u_{2m+5})F^+_n(v,w)
}.
\]
\end{rem}

\begin{rem}
In the form with apparent singularities, $B_n$ gets multiplied by
$\theta_p(xz)/\theta_p(q^n xz)$, making it elliptic, as expected.
Moreover, since for each generator of our lattice of shifts we have
exhibited an isomonodromy transformation with coef\/f\/icients of degree
independent of $n$, the same holds for an arbitrary shift.
\end{rem}

Once again, the special values, together with the determinant and the
$p$-theta law, are more than suf\/f\/icient to determine $\tilde{B}_n$; indeed,
for generic parameters, any $m+3$ of the special values suf\/f\/ice.  For both
$\tilde{A}_n$ and $\tilde{B}_n$, this gives rise to a number of relations
between the coef\/f\/icients.  For instance, using $\tilde{B}_n$, we can express
\[
G_n(u_{m+3}/q;v')F^+_n(u_{m+3}/q,w)
\]
as an explicit linear combination of the terms
\[
G^+_n(u_s/q,v')F_n(u_s;w)
\]
for $0\le s<m+3$.  If we choose $v'$ and $w$ suitably, we can arrange both
for some terms of the resulting identity to drop out, and for the remaining
integrals to be order $m$ elliptic Selberg integrals.  For instance, if
$v'=u_0/q$, $w=\widehat{u_1}$, then
\[
G^+_n(u_0/q,v') = F_n(u_1;w)=0,
\]
and thus the $s=0,1$ terms disappear, and we are left with an expression
for
\[
G_n(u_{m+3}/q;u_0/q)F_n(u_1;u_{m+3}/q)
\]
as a linear combination of
\[
G^+_n(u_s/q,u_0/q)F^-_n(u_s,u_1)
\]
for $2\le s<m+3$.  When $m=1$, these identities and the corresponding
identities arising from the 12 entry of $A_n$ give direct proofs (i.e.,
without using the $E_7$ symmetry) of new special cases of Theorem 5.1 of
\cite{recur} (which states that the elliptic Selberg integral satisf\/ies
bilinear relations making it a tau function for the elliptic Painlev\'e
equation).

The special case $\tilde{B}(z;u_6/q,u_7/q;u_0,u_1)$ (with $m=1$) is
particularly nice.  In that case, in contrast to the situation with
$\tilde{A}$, {\em all} of the integrals that appear in the expressions for
the singular values are order 1 elliptic Selberg integrals, and may thus be
expressed via \eqref{eq:eP_tau} as tau functions for elliptic Painlev\'e.
It should follow from a suitable Zariski density argument that the various
formulas resulting from consistency of this expression and of the action of
isomonodromy transformations continue to hold for arbitrary tau functions,
with any appearance of $q^n$ replaced by $Q\in \C^*$ such that
\[
Q^2 q^{-2} \prod_{0\le r<2m+6} u_r = (pq)^{m+1}.
\]
(Sketch: As $n$ varies over (large) positive integers, the balancing
condition describes a dense countable family of hypersurfaces in parameter
space; that the contour integral is dense among all solutions on such
hypersurfaces follows from the fact that its dif\/ference Galois group is
generically equal to~$\GL_2$.)  However, since the arguments of~\cite{partII} give a much more conceptual proof of this fact (and,
conversely, that any function satisfying all consistency conditions is a~tau function), there seems little point to f\/leshing out the details of a~Zariski density argument.

Also of interest in the case $m=1$ is the relation to the action of the
Weyl group~$E_7$.  This turns out to be easiest to describe in terms of the
$B$ matrices, although what meaning this has from a dif\/ference equation
perspective is as yet unclear.

\begin{thm}
In the case $m=1$, let
\[
u'_r = \begin{cases}
q^{-1/2}u_r/x, & 0\le r<4,\\
q^{1/2}u_r x, & 4\le r<8,
\end{cases}
\]
where
\[
x = \left(\frac{u_0u_1u_2u_3}{u_4u_5u_6u_7}\right)^{1/4}
  = \left(\frac{pq^{2-n}}{u_4u_5u_6u_7}\right)^{1/2}
  = \left(\frac{u_0u_1u_2u_3}{pq^{2-n}}\right)^{1/2}
.
\]
Then there exist matrices $C$ and $D$ independent of $z$ such that
\[
B_n(z;v,w;v',w';u_0,\dots,u_7)
=
C B_n\big(q^{1/2} x z;v'',w'';v''',w''';u'_0,\dots,u'_7\big) D.
\]
\end{thm}

\begin{proof}
Indeed, it suf\/f\/ices to verify this in the special case
\begin{alignat*}{5}
& v=u_6/q,\qquad &&w =u_7/q,\qquad &&v'=u_0/q, \qquad &&w'=u_1/q, &\\
& v''=u'_6/q,\qquad &&w''=u'_7/q, \qquad && v'''=u'_0/q, \qquad && w'''=u'_1/q,&
\end{alignat*}
for $z=u_r/q$, $4\le r\le 7$, in which case it follows readily by several
applications of the transformation law \cite[Theorem~9.7]{xforms}
\begin{gather*}
\II_n(u_0,\dots,u_7)
=
\prod_{1\le j\le n}\!
\left(
\prod_{0\le r<s<4} \!\Gampq(q^{n-j}u_ru_s)\!\!
\prod_{4\le r<s<8}\! \Gampq(q^{n-j}u_ru_s)\right)\!
\II_n(u'_0,\dots,u'_7).\!\!\!
\tag*{\qed}
\end{gather*}
\renewcommand{\qed}{}
\end{proof}

Finally, we observe that the case $m=0$ is precisely the elliptic
hypergeometric equation~\cite{SpiridonovVP:2005} for a (terminating)
elliptic hypergeometric series, since in that case the biorthogonal
functions are Spiridonov's elliptic hypergeometric biorthogonal functions.

\subsection*{Acknowledgements}

The author would like to thank N.~Witte for some helpful discussions of the
orthogonal polynomial approach to isomonodromy (and the University of
Melbourne for hosting the author's sabbatical when the discussions took
place), and D.~Arinkin and A.~Borodin for discussions leading to
\cite{partII} (and thus clarifying what needed (and, perhaps more
importantly, what did {\em not} need) to be established here).  The author
was supported in part by NSF grant numbered DMS-0401387, with additional
work on the project supported by NSF grants numbered DMS-0833464 and
DMS-1001645.

\pdfbookmark[1]{References}{ref}
\LastPageEnding


\begin{thebibliography}{99}

\footnotesize\itemsep=0pt

\bibitem{AdlerM/vanMoerbekeP:1999}
Adler M., van Moerbeke P.,
The spectrum of coupled random matrices,
\href{http://dx.doi.org/10.2307/121077}{{\em Ann. of Math.~(2)}} {\bf 149} (1999), 921--976,
\href{http://arxiv.org/abs/hep-th/9907213}{hep-th/9907213}.

\bibitem{ArinkinD/BorodinA:2006}
Arinkin D., Borodin A.,
Moduli spaces of $d$-connections and dif\/ference Painlev\'e   equations,
\href{http://dx.doi.org/10.1215/S0012-7094-06-13433-6}{{\em Duke Math.~J.}} {\bf 134} (2006), 515--556,
\href{http://arxiv.org/abs/math.AG/0411584}{math.AG/0411584}.

\bibitem{partII}
Arinkin~D., Borodin~A.,   Rains E.M.,
in preparation.

\bibitem{BirkhoffGD:1911}
 Birkhof\/f G.D.,
General theory of linear dif\/ference equations,
\href{http://dx.doi.org/10.2307/1988577}{{\em Trans. Amer. Math. Soc.}} {\bf 12} (1911), 243--284.

\bibitem{BorodinA:2004}
Borodin A.,
Isomonodromy transformations of linear systems of dif\/ference   equations,
\href{http://dx.doi.org/10.4007/annals.2004.160.1141}{{\em Ann. of Math.~(2)}} {\bf 160} (2004), 1141--1182,
\href{http://arxiv.org/abs/math.CA/0209144}{math.CA/0209144}.

\bibitem{DeiftPA:1999}
 Deift P.A.,
Orthogonal polynomials and random matrices: a Riemann--Hilbert approach,   {\em Courant Lecture Notes in
  Mathematics}, Vol.~3,  New York University, Courant Institute of Mathematical Sciences, New   York, 1999.

\bibitem{vanDiejenJF/SpiridonovVP:2000}
  van~Diejen J.F., Spiridonov V.P.,
An elliptic Macdonald--Morris conjecture and multiple modular hypergeometric sums,
  {\em Math. Res. Lett.} {\bf 7} (2000), 729--746.

\bibitem{EtingofPI:1993}
 Etingof P.I.,
 Dif\/ference equations with elliptic coef\/f\/icients and quantum af\/f\/ine   algebras,
\href{http://arxiv.org/abs/hep-th/9312057}{hep-th/9312057}.

\bibitem{EtingofPI:1995}
 Etingof P.I.,
Galois groups and connection matrices of $q$-dif\/ference equations,
\href{http://dx.doi.org/10.1090/S1079-6762-95-01001-8}{{\em Electron. Res. Announc. Amer. Math. Soc.}} {\bf 1} (1995), no.~1, 1--9.


\bibitem{FokasAS/ItsAR/KitaevAV:1992}
 Fokas A.S., Its A.R., Kitaev A.V.,
The isomonodromy approach to matrix models in 2D quantum gravity,
\href{http://dx.doi.org/10.1007/BF02096594}{{\em Comm. Math. Phys.}} {\bf 147} (1992), 395--430.

\bibitem{ForresterPJ/WitteNS:2006}
 Forrester P.J., Witte N.S.,
Bi-orthogonal polynomials on the unit circle, regular semi-classical weights and integrable systems,
\href{http://dx.doi.org/10.1007/s00365-005-0616-7}{{\em Constr. Approx.}} {\bf 24} (2006), 201--237,
\href{http://arxiv.org/abs/math.CA/0412394}{math.CA/0412394}.

\bibitem{KajiwaraK/MasudaT/NoumiM/OhtaY/YamadaY:2003}
Kajiwara K., Masuda T., Noumi M., Ohta Y., Yamada Y.,
${}\sb {10}E\sb 9$ solution to the elliptic Painlev\'e equation,
\href{http://dx.doi.org/10.1088/0305-4470/36/17/102}{{\em J.~Phys.~A: Math. Gen.}} {\bf 36} (2003), L263--L272,
\href{http://arxiv.org/abs/nlin.SI/0303032}{nlin.SI/0303032}.

\bibitem{KitaevAV:2000}
 Kitaev A.V.,
Special functions of the isomonodromy type,
\href{http://dx.doi.org/10.1023/A:1006390032014}{{\em Acta Appl. Math.}} {\bf 64} (2000), 1--32.

\bibitem{KricheverIM:2004}
 Krichever I.M.,
Analytic theory of dif\/ference equations with rational and elliptic coef\/f\/icients and the Riemann--Hilbert problem,
{\it Uspekhi Mat. Nauk} {\bf 59} (2004), no.~6, 111--150 (English transl.: \href{http://dx.doi.org/10.1070/RM2004v059n06ABEH000798}{{\it Russian Math. Surveys}} {\bf 59} (2004), no.~6, 1117--1154),
\href{http://arxiv.org/abs/math-ph/0407018}{math-ph/0407018}.


\bibitem{MagnusAP:1995}
 Magnus A.P.,
Painlev\'e-type dif\/ferential equations for the recurrence coef\/f\/icients of semi-classical orthogonal polynomials,
\href{http://dx.doi.org/10.1016/0377-0427(93)E0247-J}{{\em J. Comput. Appl. Math.}} {\bf 57} (1995), 215--237,
\href{http://arxiv.org/abs/math.CA/9307218}{math.CA/9307218}.

\bibitem{PraagmanC:1986}
Praagman C.,
Fundamental solutions for meromorphic linear dif\/ference equations in the complex plane, and related problems,
\href{http://dx.doi.org/10.1515/crll.1986.369.101}{{\em J. Reine Angew. Math.}} {\bf 369} (1986), 101--109.

\bibitem{vanderPutM/SingerMF:1997}
van~der Put M.,  Singer M.F.,
Galois theory of dif\/ference equations,
{\em Lecture Notes in Mathematics}, Vol.~1666, Springer-Verlag, Berlin, 1997.

\bibitem{recur}
Rains E.M.,
Recurrences of elliptic hypergeometric integrals,
in  Elliptic Integrable Systems, Editors M.~Noumi and K.~Takasaki, {\em Rokko Lectures in Mathematics}, Vol.~18, Kobe University, 2005, 183--199,
\href{http://arxiv.org/abs/math.CA/0504285}{math.CA/0504285}.

\bibitem{xforms}
Rains E.M.,
Transformations of elliptic hypergeometric integrals,
\href{http://dx.doi.org/10.4007/annals.2010.171.169}{{\em Ann. of Math.~(2)}} {\bf 171} (2010), 169--243,
\href{http://arxiv.org/abs/math.QA/0309252}{math.QA/0309252}.

\bibitem{RuijsenaarsSNM:1999}
 Ruijsenaars S.N.M.,
A generalized hypergeometric function satisfying four analytic dif\/ference equations of Askey--Wilson type,
\href{http://dx.doi.org/10.1007/s002200050840}{{\em Comm. Math. Phys.}} {\bf 206} (1999), 639--690.

\bibitem{SakaiH:2001}
Sakai H.,
Rational surfaces associated with af\/f\/ine root systems and geometry of the Painlev\'e equations,
\href{http://dx.doi.org/10.1007/s002200100446}{{\em Comm. Math. Phys.}} {\bf 220} (2001), 165--229.

\bibitem{SpiridonovVP:2001b}
 Spiridonov V.P.,
Elliptic beta integrals and special functions of hypergeometric type,
in  Integrable Structures of Exactly Solvable Two-Dimensional  Models of Quantum Field Theory (Kiev, 2000),
{\em NATO Sci. Ser. II Math. Phys. Chem.}, Vol.~35, Kluwer Acad. Publ., Dordrecht, 2001, 305--313.

\bibitem{SpiridonovVP:2001a}
Spiridonov V.P.,
On the elliptic beta function,
{\em Uspekhi Mat. Nauk} {\bf 56} (2001), no.~1, 181--182
(English transl.: \href{http://dx.doi.org/10.1070/rm2001v056n01ABEH000374}{{\it Russian Math. Surveys}} {\bf 56} (2001), no.~1, 185--186).

\bibitem{SpiridonovVP:2005}
Spiridonov V.P.,
Classical elliptic hypergeometric functions and their applications,
in  Elliptic Integrable Systems, Editors M.~Noumi and K.~Takasaki, {\em Rokko Lectures in Mathematics}, Vol.~18, Kobe University, 2005, 253--288,
\href{http://arxiv.org/abs/math.CA/0511579}{math.CA/0511579}.

\bibitem{SzegoG:1939}
Szeg\"o G.,
Orthogonal polynomials,
{\it American Mathematical Society Colloquium Publications}, Vol.~23, Ame\-ri\-can Mathematical Society, New York, 1939.


\bibitem{YamadaY:2009}
Yamada Y.,
A Lax formalism for the elliptic dif\/ference Painlev\'e equation,
\href{http://dx.doi.org/10.3842/SIGMA.2009.042}{{\em SIGMA}} {\bf 5} (2009), 042,   15~pages,
\href{http://arxiv.org/abs/0811.1796}{arXiv:0811.1796}.

  \end{thebibliography}
\end{document}